Yu. V. Troshchiev

Improvement of the monotonicity properties of the difference schemes by building in them of the monotonizing operators

M.V. Lomonosov Moscow State University, Department of Calculate Mathematics and Cybernetics.

E-mail: yuvt@yandex.ru.

**Abstract.**

The method of monotonization of difference schemes is being considered in the paper. The method was earlier proposed by the author for stationary problems. It is investigated in the paper more profoundly. The idea of the method is to build the monotonizing operators into the schemes so that the balance relations from point to point are not violated. Different monotonizing operators can be used to be installed in the schemes. Propositions concerning approximation and stability of the monotonized schemes are formulated and proved. Also a proposition significant for practical use of the schemes is formulated and proved. The idea is to use the monotonized schemes in the cases when the proposition conditions are fulfilled. The proposition is based on closeness of solutions of the initial and auxiliary schemes. Constructions for solving of time dependent problems are also written in the paper. One dimensional example and three-dimensional hydrodynamic example are considered. The method allows to considerably decrease value of calculations in many cases.

## 1. Method construction.

Let

$$\mathbf{A}(\mathbf{u}) = 0 \tag{1}$$

is some difference scheme [1-4] for a stationary boundary-value problem. $\mathbf{u}$ – is a mesh function

$$\mathbf{u} \in R^n, \tag{2}$$

$R^n$ – real $n$-dimensional space, $\mathbf{A}$ is a finite dimensional operator (may be nonlinear):

$$\mathbf{A}: R^n \to R^n. \tag{3}$$

It is supposed that the problem may contain several variables and be multidimensional. For example if $\mathbf{u}$ describes a three dimensional velocity $\mathbf{v} = (v_x, v_y, v_z)$ and a pressure $p$ in a three-dimensional parallelepiped then the dimension $n$ equals $4m$:

$$\mathbf{u} = (v_{x1},...,v_{xm}, v_{y1},...,v_{ym}, v_{z1},...,v_{zm}, p_1,...,p_m)^T. \tag{4}$$

Here $m = N_x N_y N_z$, where $N_x$, $N_y$ and $N_z$ are mesh dimensions in the correspondent directions. Such notification is given here to represent a physical problem in the form (1). Further we shall use these single-indexed notification or multi-indexed notification of the variables of the multi-dimensional problems depending on the context and being ensured that a one-to-one correspondence is established between these two forms. It is supposed here that the scheme is constructed by the balance method or by the finite-difference approximation method.

The solutions of such problems are often nonmonotonic (oscillating from point to point), and the mesh should be refined more than it is necessary to describe the being investigated function adequately. Note, that a nonmonotonicity means a nonmonotonicity in space (not in time) here, and it is the nonmonotonicity from point to point, i.e. oscillations take place at each or at almost each mesh step. And the exact continuous solution, for example in one dimensional case, may be nonmonotonic but it is monotonic in the intervals between extremes.

Methods of construction of monotonic difference schemes, different variants of smoothing operators for monotonization of the found solutions and some other related topics are considered in the papers [3-15] and others. Often such a monotonization disturbs balance equations. Some method of monotonization conserving balance equations was proposed by the author [16, 17]. This method is being investigated in this paper. The idea of the method is not to monotonize the found solution but to install the monotonizing operator into the difference scheme so that the solution will be already monotonized. Practically it is often possible to get similar results by monotonization of the found solution [9] despite the disturbance of the balance equations. Another way is to construct initially monotonous schemes. Comparison with the Flux-Corrected Methods [3,8] is more difficult because of their variety. May be some of such methods can be transformed to the form that discussed in this paper. But the method proposed in [16,17] is based on other ideas. It is mathematical and doesn't use physical considerations.

Let us describe the construction of the monotonization method. Note that here we speak not only about an unconditional monotonization but also about improvement of monotonicity properties in comparison with the initial scheme. The difference scheme (1) can be written in the form

$$\mathbf{F}(\mathbf{D}(\mathbf{u})) = 0 \qquad (5)$$

where the operator $\mathbf{D}$ calculates difference derivatives up to the necessary order and the operator $\mathbf{F}$ describes physical laws. Usually operator $\mathbf{D}$ is «explicit» – it calculates values of the derivatives by explicit formulas using several surrounding points of some pattern. If the solution of the scheme (5) is nonmonotonic then it is offered [16,17] to improve monotonicity properties by changing the operator $\mathbf{D}$ by some "implicit" operator $\mathbf{D}_{imp}$

$$\mathbf{F}(\mathbf{D}_{imp}(\mathbf{u})) = 0, \qquad (6)$$

which calculates the vector of derivatives from the vector of the mesh function as a solution of equation, i.e. there is an analogy with implicit schemes on time. Note here S.K. Godunov scheme

"cancellation of discontinuity" and TVD-scheme of Chakravarthy-Osher that produce monotonic solutions and calculate derivatives in different ways depending on some conditions [3,5,6,12].

Improvement of monotonicity properties can be considered as disappearing of the false extreme points, decreasing of oscillations amplitude, decreasing of the total fluctuation [12]. The operator $\mathbf{D}_{imp}$ can be constructed of the operator $\mathbf{D}$ and some monotonizing operator $\mathbf{M}$ applying some algorithm which will be described in what fallows. As a consequence of this procedure the monotonicity properties of the solution become improved under some usually fulfilled conditions.

"Explicit" operators $\mathbf{B}$ and $\mathbf{C}$ can be constructed so, that

$$\mathbf{F}(\mathbf{Bv}) = 0$$
$$\mathbf{u} = \mathbf{Cv}. \tag{7}$$

It leads to decrease of calculations amount. Here $\mathbf{v}$ is an auxiliary mesh function, and $\mathbf{u}$ is the desired solution.

Constructively the operator $\mathbf{D}_{imp}$ is derived in the following manner. Let for example a second order differential equation is under consideration. Let us write the scheme (5) in the form

$$\mathbf{F}(\mathbf{u}, D_1\mathbf{u}, D_2\mathbf{u}) = 0, \tag{8}$$

where $\mathbf{D}_1$ and $\mathbf{D}_2$ are the operators that calculate the first and the second derivatives, correspondingly.

Monotonizing (smoothing) operator in three-dimensional case like (4) with a regular mesh at each direction can be for example of the following form:

$$(\mathbf{Mu})_{i,j,k} = \frac{1}{6}\left( \frac{u_{i+1,j,k} + u_{i,j,k}}{2} + \frac{u_{i-1,j,k} + u_{i,j,k}}{2} + \right.$$
$$\frac{u_{i,j+1,k} + u_{i,j,k}}{2} + \frac{u_{i,j-1,k} + u_{i,j,k}}{2} +$$
$$\left. \frac{u_{i,j,k+1} + u_{i,j,k}}{2} + \frac{u_{i,j,k-1} + u_{i,j,k}}{2} \right) . \tag{9}$$

Denotation $(\mathbf{Mu})_{i,j,k}$ means the component $w_{i,j,k}$ of the vector $\mathbf{w} = \mathbf{Mu}$, and the symbol $u$ takes values $v_x$, $v_y$, $v_z$ and $p$. Some complexity of the denotation caused by the fact that in three-dimensional case it is natural to understand vector $\mathbf{u}$ as a vector consisting of vectors ($v_x, v_y, v_z, p$). Besides, the ranges of indexes alteration are not pointed out in (9) because they depend on mesh and region shape under consideration. We suppose that formula is defined in all those points $(i, j, k)$ in which all the variables incorporated in the formula are defined. Another monotonizing operators can also be used. Two properties of the operator $\mathbf{M}$ are essential in this paper: 1) if $\mathbf{v}$ is a nonmonotonic mesh function, then monotonicity properties of the function $\mathbf{Mv}$ are better, than monotonicity properties of the function $\mathbf{v}$; 2) if functions $\mathbf{u}$ and $\mathbf{v}$ differ little then monotonicity properties of the function $\mathbf{Mv}$ are also better than monotonicity properties of the function $\mathbf{u}$. Smoothing and averaging operators that reduce oscillations of mesh functions are being discussed for example in [7,9,10]. Variants of definition of monotonicity properties are discussed for example in [12,17].

To construct the monotonized scheme let us change the first argument in the scheme (8) by means of operator (9):

$$\mathbf{F}(\mathbf{Mv}, \mathbf{D}_1\mathbf{v}, \mathbf{D}_2\mathbf{v}) = 0. \qquad (10)$$

As a result an auxiliary difference scheme is derived in which not only derivatives are approximated but the function is approximated as well. To distinguish the solution of this scheme from the solution $\mathbf{u}$ of the scheme (8) we denote it $\mathbf{v}$. But balance equations are written relatively to the function $\mathbf{Mv}$ now. So we shall consider this function $\mathbf{Mv}$ as a solution

$$\mathbf{y} = \mathbf{Mv}. \qquad (11)$$

It is convenient to solve the scheme (10), (11) in such a form but it is possible to write it in the usual form similar to (8) if $\mathbf{M}$ is invertible:

$$\mathbf{F}(\mathbf{y}, \mathbf{D}_1\mathbf{M}^{-1}\mathbf{y}, \mathbf{D}_2\mathbf{M}^{-1}\mathbf{y}) = 0. \qquad (12)$$

It is necessary to find derivatives as solutions of equations in such a form but may be it is more convenient to derive proofs for such a notation.

Schemes (10), (11) and (12) are a new kind of schemes. In (10), (11) the operator $\mathbf{M}$ can be degenerate and in (12) it should be nondegenerate. So the scheme (10), (11) is more widely applicable.

*Proposition 1.* Approximation order of the scheme (10), (11) is greater or equal than the smaller of two numbers – approximation order of the scheme (8) and approximation order of monotonizing operator.

Proof is evident.

*Proposition 2.* Let the scheme (8) is linear. If there exists such $h_0 > 0$ that the linear operators $(\partial \mathbf{F}/\partial \mathbf{D}_2)\mathbf{D}_2$ and $\mathbf{M}$ are nondegenerate for any $h < h_0$ then the scheme (12) is stable. $h$ here is $\max(h_x, h_y, h_z)$, $h_x$, $h_y$, $h_z$ – mesh steps at $x$, $y$ and $z$ directions, correspondingly. $(\partial \mathbf{F}/\partial \mathbf{D}_2)\mathbf{D}_2$ is a linearization of $\mathbf{F}$ respectively to the third argument.

*Proof.* Since the operators $(\partial \mathbf{F}/\partial \mathbf{D}_2)\mathbf{D}_2$ and $\mathbf{M}$ are nondegenerate for any $h < h_0$ the operator $(\partial \mathbf{F}/\partial \mathbf{D}_2)\mathbf{D}_2\mathbf{M}^{-1}$ is also nondegenerate for these $h$. The left side of (8) is $\partial \mathbf{F}/\partial \mathbf{u} + (\partial \mathbf{F}/\partial \mathbf{D}_1)\mathbf{D}_1 + (\partial \mathbf{F}/\partial \mathbf{D}_2)\mathbf{D}_2$ and the left side of (12) is $\partial \mathbf{F}/\partial \mathbf{u} + (\partial \mathbf{F}/\partial \mathbf{D}_1)\mathbf{D}_1\mathbf{M}^{-1} + (\partial \mathbf{F}/\partial \mathbf{D}_2)\mathbf{D}_2\mathbf{M}^{-1}$. The first member $\partial \mathbf{F}/\partial \mathbf{u}$ in both operators doesn't depend on $h$, so $\|\partial \mathbf{F}/\partial \mathbf{u}\| = C_1 = \text{const}$. The second members $(\partial \mathbf{F}/\partial \mathbf{D}_1)\mathbf{D}_1$ and $(\partial \mathbf{F}/\partial \mathbf{D}_1)\mathbf{D}_1\mathbf{M}^{-1}$ are proportional to $h^{-1}$: $\|(\partial \mathbf{F}/\partial \mathbf{D}_1)\mathbf{D}_1\| = C_2 h^{-1}$ and $\|(\partial \mathbf{F}/\partial \mathbf{D}_1)\mathbf{D}_1\mathbf{M}^{-1}\| = \widetilde{C}_2 h^{-1}$. The third members $(\partial \mathbf{F}/\partial \mathbf{D}_2)\mathbf{D}_2$ and $(\partial \mathbf{F}/\partial \mathbf{D}_2)\mathbf{D}_2\mathbf{M}^{-1}$ are nondegenerate and proportional to $h^{-2}$: $\|(\partial \mathbf{F}/\partial \mathbf{D}_2)\mathbf{D}_2\| = C_3 h^{-2}$ and $\|(\partial \mathbf{F}/\partial \mathbf{D}_2)\mathbf{D}_2\mathbf{M}^{-1}\| = \widetilde{C}_3 h^{-2}$. So (8) and (12) are nondegenerate for small enough $h$ and the schemes are stable. The proof is over.

*Note.* Matters of approximation and convergence are little incorrect in this context because approximation and convergence are limit concepts. But we here discuss a possibility of calculations with a large enough step, i.e. without converging of the step value to zero.

Let now some positive functional $f(\mathbf{u})$ is defined at the set of mesh functions $\{\mathbf{u}\}$ and this functional numerically characterizes the value of nonmonotonicity of the functions. Let this functional is Lipschitzian in $\|\cdot\|_C$ norm:

$$f(\mathbf{u}): \quad |f(\mathbf{u}) - f(\mathbf{v})| \leq K \|\mathbf{u} - \mathbf{v}\|_C. \tag{13}$$

*Proposition 3.* Let at some $h$ $f(\mathbf{u}) = \delta > 0$ where $\mathbf{u}$ is a solution of the scheme (8) and $f(\mathbf{Mu}) = k\delta$, $0 < k < 1$. Let also $\|\mathbf{u} - \mathbf{v}\|_C = \varepsilon > 0$, $\varepsilon < \delta$, $K \|\mathbf{M}\|_C \varepsilon < k\delta$ at this $h$ where $\mathbf{v}$ is a solution of (10). Then $f(\mathbf{Mv}) = k_1 \delta$ where

$$\frac{k\delta - K \|\mathbf{M}\|_C \varepsilon}{\delta + \varepsilon} < k_1 < \frac{k\delta + K \|\mathbf{M}\|_C \varepsilon}{\delta - \varepsilon}. \tag{14}$$

*Proof.* $\|\mathbf{Mu} - \mathbf{Mv}\|_C \leq \|\mathbf{M}\|_C \cdot \|\mathbf{u} - \mathbf{v}\|_C = \|\mathbf{M}\|_C \varepsilon$. Then $|f(\mathbf{Mu}) - f(\mathbf{Mv})| \leq K \|\mathbf{Mu} - \mathbf{Mv}\|_C \leq K \|\mathbf{M}\|_C \varepsilon$. So

$$f(\mathbf{v}) \in (\delta - \varepsilon, \delta + \varepsilon),$$

$$f(\mathbf{Mv}) \in f(\mathbf{Mv}) \in (k\delta - K \|\mathbf{M}\|_C \varepsilon, k\delta + K \|\mathbf{M}\|_C \varepsilon)$$

and $f(\mathbf{Mv}) = k_1 \delta$ where $k_1$ is from the necessary interval. The proof is over.

*Note 1.* Let the scheme (8) is linear and there exists such $h_0 > 0$ that the linear operator $(\partial \mathbf{F} / \partial \mathbf{D}_2) \mathbf{D}_2$ is nondegenerate for any $h_1 < h_0$. Then $\varepsilon$ converges to zero like $h^2$.

*Note 2.* The norm of monotonizing operators is usually of unit order. If functional $f$ means maximal change per mesh step then $K = 2$. $\delta$ for an oscillating function can be for example greater than $h$ and $\varepsilon$ in accordance with a scheme precision can be like $h^2$. So satisfaction of Proposition 3 conditions is a usual situation. The idea is to use monotonized schemes in cases when the proposition conditions are fulfilled.

It may be more difficult to distinguish operators $\mathbf{D}_1$ and $\mathbf{D}_2$ in nonlinear schemes and schemes constructed by balance method. But there also exist terms proportional to $1/h$ and $1/h^2$. And some attention to physical sense is need in these cases. Later there is an example of monotonization of nonlinear scheme.

The possibility to find monotonic solution at less number of mesh points is essential in multidimensional problems. For example decreasing of the number of points by a factor of 5 in each direction for a three-dimensional problem causes decreasing of the total number of points by a factor of 125. Practically the necessary number of points usually can be determined by small changing of solution under reducing of the mesh step by factor of 2. Note that schemes considered in the paper [9] algebraically are not (10), (11) but (8) and further smoothing of the solution.

Finally in this item let us construct conservative monotonized scheme for dynamical difference scheme

$$\frac{\mathbf{u}^{n+1} - \mathbf{u}^n}{\tau} = \sigma \mathbf{F}(\mathbf{u}^{n+1}, \mathbf{D}_1\mathbf{u}^{n+1}, \mathbf{D}_2\mathbf{u}^{n+1}) + (1-\sigma)\mathbf{F}(\mathbf{u}^n, \mathbf{D}_1\mathbf{u}^n, \mathbf{D}_2\mathbf{u}^n), \quad \sigma \in [0,1]. \quad (15)$$

Applying the form (12) we derive

$$\begin{aligned}\frac{\mathbf{y}^{n+1} - \mathbf{y}^n}{\tau} &= \\ &= \sigma \mathbf{F}(\mathbf{y}^{n+1}, \mathbf{D}_1\mathbf{M}^{-1}\mathbf{y}^{n+1}, \mathbf{D}_2\mathbf{M}^{-1}\mathbf{y}^{n+1}) + (1-\sigma)\mathbf{F}(\mathbf{y}^n, \mathbf{D}_1\mathbf{M}^{-1}\mathbf{y}^n, \mathbf{D}_2\mathbf{M}^{-1}\mathbf{y}^n)\end{aligned} \quad (16)$$

and using the auxiliary variable $\mathbf{v}$ we derive

$$\begin{aligned}\frac{\mathbf{y}^{n+1} - \mathbf{M}\mathbf{v}^n}{\tau} &= \sigma \mathbf{F}(\mathbf{M}\mathbf{v}^{n+1}, \mathbf{D}_1\mathbf{v}^{n+1}, \mathbf{D}_2\mathbf{v}^{n+1}) + (1-\sigma)\mathbf{F}(\mathbf{M}\mathbf{v}^n, \mathbf{D}_1\mathbf{v}^n, \mathbf{D}_2\mathbf{v}^n), \\ \mathbf{v}^{n+1} &= \mathbf{M}^{-1}\mathbf{y}^{n+1}.\end{aligned} \quad (17)$$

or in another form

$$\begin{aligned}\mathbf{v}^{n+1} &= \mathbf{v}^n + \tau\mathbf{M}^{-1}\sigma\mathbf{F}(\mathbf{M}\mathbf{v}^{n+1}, \mathbf{D}_1\mathbf{v}^{n+1}, \mathbf{D}_2\mathbf{v}^{n+1}) + (1-\sigma)\mathbf{F}(\mathbf{M}\mathbf{v}^n, \mathbf{D}_1\mathbf{v}^n, \mathbf{D}_2\mathbf{v}^n), \\ \mathbf{y}^{n+1} &= \mathbf{M}\mathbf{v}^{n+1}.\end{aligned} \quad (18)$$

In all the forms (16)-(18) if $\sigma \neq 0$ it is necessary to solve an equation

$$\mathbf{M}\mathbf{a} = \mathbf{b} \tag{19}$$

at each iteration to find $\mathbf{y}^{n+1}$ or $\mathbf{v}^{n+1}$.

**2. One-dimensional scheme.**

Let us consider a boundary value problem:

$$\begin{aligned} &k_0 + k_1 U + k_2 U' + k_3 U'' = 0, \\ &U(a) = u_0, \quad U(b) = u_{n+1}. \end{aligned} \tag{18}$$

Now we construct a difference scheme for the problem (18). Let $\mathbf{u}$ is a mesh function $\{u_i\}$, $i = 1,2,...,n$ defined at the regular mesh $\{x_i\}$, $i = 1,2,...,n$, $x_{i+1} - x_i = h$, $i = 0,2,...,n$. The values $u_0$ and $u_{n+1}$ we consider as given boundary conditions. Let operators for calculation of the first and the second derivatives be the following:

$$\begin{aligned} (\mathbf{D}_1 \mathbf{u})_i &= \frac{1}{h}(\widetilde{\mathbf{D}}_1 \mathbf{u})_i = (u_{i+1} - u_{i-1})/2h, \\ (\mathbf{D}_2 \mathbf{u})_i &= \frac{1}{h^2}(\widetilde{\mathbf{D}}_2 \mathbf{u})_i = (u_{i+1} - 2u_i + u_{i-1})/h^2, \quad i = 1,2,...,n, \end{aligned} \tag{19}$$

so we write the following difference scheme:

$$h^2 k_0 + h^2 k_1 u_i + h k_2 (\widetilde{\mathbf{D}}_1 \mathbf{u})_i + k_3 (\widetilde{\mathbf{D}}_2 \mathbf{u})_i = 0, \quad i = 1,2,...,n, \tag{20}$$

where $k_0$, $k_1$, $k_2$, $k_3$ are real numbers. For the problem to be nondegenerate we assume $k_3 \neq 0$. Thus a difference scheme of the second order approximation is constructed for the problem (18). Let us improve monotonicity properties of this scheme.

*Definition.* We shall speak that a mesh function $\mathbf{u}$ oscillates from point to point at the interval from $k$ to $l$ ($k \leq i \leq l$) when and only when inequalities $u_{i-1} > u_i$ and $u_{i+1} > u_i$ are satisfied for all even (odd) $i$ from the interval $k+1 \leq i \leq l-1$ and inequalities $u_{i-1} < u_i$ and $u_{i+1} < u_i$ are satisfied for all odd (even) $i$ from this interval.

The monotonizing operator $\mathbf{M}$ let be the following

$$(\mathbf{M}\mathbf{u})_i = ((u_{i+1} + u_i)/2 + (u_i + u_{i-1})/2)/2 = (u_{i+1} + 2u_i + u_{i-1})/4. \tag{21}$$

Note that such an operator is brought in the text-book [10] to filter pressure oscillations for calculation of their amplitude.

Let the improvement of monotonicity properties be for example a decrease of maximal change of mesh function at the mesh steps:

$$\max_{i=k,k+1,\ldots,I-1} \left(|(\mathbf{Mu})_{i+1} - (\mathbf{Mu})_i|\right) < \max_{i=k,k+1,\ldots,I-1} \left(|u_{i+1} - u_i|\right). \tag{22}$$

*Note 1.* The function $\mathbf{Mu}$ may be nonoscillatory at the considered interval where the function $\mathbf{u}$ oscillates.

*Note 2.* Precise clarification of the operator $\mathbf{M}$ properties is not necessary in this paper but usually the fact is that it improves monotonicity properties of the oscillating from point to point functions.

Now we write the auxiliary difference scheme:

$$h^2 k_0 + h^2 k_1 (\mathbf{Mv})_i + h k_2 (\widetilde{\mathbf{D}}_1 \mathbf{v})_i + k_3 (\widetilde{\mathbf{D}}_2 \mathbf{v})_i = 0, \quad i = 1,2,\ldots,n, \tag{23}$$

and the monotonized solution

$$y_i = (\mathbf{Mv})_i \quad i = 1,2,\ldots,n. \tag{24}$$

It is possible to write difference scheme (20) in the form (8)

$$\mathbf{F}(\mathbf{u}, \mathbf{D}_1 \mathbf{u}, \mathbf{D}_2 \mathbf{u}) = 0 \tag{25}$$

And the difference scheme (23), (24) in the form (10), (11)

$$\mathbf{F}(\mathbf{Mv}, \mathbf{D}_1 \mathbf{v}, \mathbf{D}_2 \mathbf{v}) = 0,$$
$$\mathbf{y} = \mathbf{Mv}. \tag{26}$$

Solutions of the systems (23) and (20) differ the following

$$h^2 k_0 \mathbf{E}(\mathbf{u} - \mathbf{v}) + h^2 k_1 (\mathbf{Eu} - \mathbf{Mv}) + \left(h k_2 \widetilde{\mathbf{D}}_1 + k_3 \widetilde{\mathbf{D}}_2\right)(\mathbf{u} - \mathbf{v}) = 0. \tag{27}$$

So if $\mathbf{u}$ is nonmonotonic then conditions analogous to proposition 3 ones are also normally satisfied in this case.

*Note 1.* Verification of the conditions can be done for example numerically.

*Note 2.* Solution of the scheme (23), (24) will be essentially more precise in the case if oscillations take place around the precise solution.

An example of such a situation for a mesh of 11 points is given at fig. 1. Comparing with [17] $k_3$ is not small here against $k_1$.

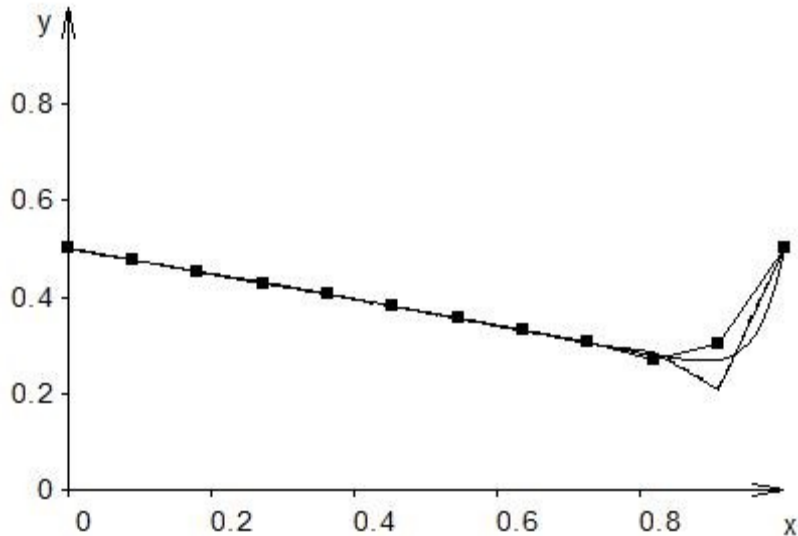

**Fig. 1** Solutions of the eq. (20) – solid polygonal line, eq. (23) –it is almost coincide with the previous line and so invisible, and eq. (23), (24) – marked line at parameters values: $x_1=0$, $x_4=1$, $k_0=10$, $k_1=-5$, $k_2=30$, $k_3=-1$, $u_1=0.5$, $u_4=0.5$. Solid curve is a solution at 100 point mesh which can be considered as almost precise one

In [17] it is sown that smallness of $k_1 h^2$ is not a necessary condition of more monotonicity of the function **y** in comparison with the function **u**. Also it is shown there that determinant of (23) can be zero for some $h$ when determinant of (20) is nonzero, i.e. the solution of (23), (24) is not always better than the solution of (20).

### 3. Difference scheme for steady equations of incompressible liquid flow.

It is possible to make analogous reasoning like in previous paragraph for linear multidimensional case.

Check-up of the construction (10), (11) for steady equations of incompressible liquid flow showed that monotonicity properties of the solution essentially improve. Calculations were performed in a cubic area at a regular mesh but the method can be generalized for many other cases. To illustrate a calculation of a nonlinear physical problem using the scheme (10), (11) let us consider a problem that was investigated by the author while investigation of liquid flows through filter cells [18,19]. Consider a system of Navier-Stokes equations:

$$v_x \frac{\partial v_x}{\partial x} + v_y \frac{\partial v_x}{\partial y} + v_z \frac{\partial v_x}{\partial z} = -\frac{1}{\rho}\frac{\partial p}{\partial x} + \nu \Delta v_x,$$
$$v_x \frac{\partial v_y}{\partial x} + v_y \frac{\partial v_y}{\partial y} + v_z \frac{\partial v_y}{\partial z} = -\frac{1}{\rho}\frac{\partial p}{\partial y} + \nu \Delta v_y,$$
$$v_x \frac{\partial v_z}{\partial x} + v_y \frac{\partial v_z}{\partial y} + v_z \frac{\partial v_z}{\partial z} = -\frac{1}{\rho}\frac{\partial p}{\partial z} + \nu \Delta v_z,$$
$$\frac{\partial^2 v_x}{\partial x^2} + \frac{\partial^2 v_y}{\partial y^2} + \frac{\partial^2 v_z}{\partial z^2} = 0.$$

(28)

Here $v(x,y,z,t) = (v_x, v_y, v_z)$ is liquid velocity, $p$ – pressure, $\rho$ – density, $\nu$ – viscosity. A solution in a cubic area $x \in [0,L]$, $y \in [0,L]$, $z \in [0,L]$ is being considered. There are equal square holes on the opposite facets of the cube parallel to the plane $Oyz$ and pressure functions are given in the sections of these holes

$$p(0,y,z) = p_0, \ p(L,y,z) = p_1, \ x,y \in [\frac{L}{2}-a, \frac{L}{2}+a]. \tag{29}$$

Zero velocities are established at cube facets excepting the holes:

$$v_x = v_y = v_z = 0. \tag{30}$$

And zero velocity derivatives are established at holes sections:

$$\frac{\partial v_x}{\partial x} = \frac{\partial v_y}{\partial x} = \frac{\partial v_z}{\partial x} = 0. \tag{31}$$

Such boundary conditions are physically possible but they can cause nonmonotonicity of the solution under a large mesh step.

Let us construct a simple difference scheme of the first order of approximation to avoid a detailed description. Let $N$ is some natural number which defines a regular mesh with a step $h$ for all three directions:

$$h = L/N: \; x_i = ih, \; y_j = jh, \; z_k = kh, \; i,j,k = 0,1,\ldots,N, \tag{32}$$

The mesh decomposes the area into $N^3$ cubic cells. Let us search velocities and pressures in the centers of the cells. Approximation of the first and the second derivatives were similar to (19) for each direction. To approximate the first and the second derivatives near the boundaries the boundary conditions for the functions or their derivatives on the cell facets were used. To approximate pressure derivatives near the walls at orthogonal directions inner differences (right or left ones) were used.

To solve the equations (28) in the difference form the following iterative procedure was performed:

$$\mathbf{v}^{n+1} = \mathbf{v}^n + \sigma_v \left( -v_x^n \frac{\partial \mathbf{v}^n}{\partial x} - v_y^n \frac{\partial \mathbf{v}^n}{\partial y} - v_z^n \frac{\partial \mathbf{v}^n}{\partial z} - \frac{1}{\rho} \frac{\partial p^n}{\partial x} + \nu \Delta \mathbf{v}^n \right),$$
$$p^{n+1} = p^n + \sigma_p \mathrm{div}\, \mathbf{v}^{n+1}. \tag{33}$$

Here the variables are considered as net variables $\mathbf{v} = (v_x, v_y, v_z) = \{\mathbf{v}_{ijk}\}$, $p = \{p_{ijk}\}$ and derivatives are considered as difference derivatives as described above. $\sigma_v$ and $\sigma_p$ are some iterative parameters, $n$ − number of iteration. After this procedure converged the variable $p$ can be considered as a dependent variable so we shall monotonize only $\mathbf{v}$. The monotonized scheme is the following

$$\mathbf{M}v_x \frac{\partial \mathbf{v}}{\partial x} + \mathbf{M}v_y \frac{\partial \mathbf{v}}{\partial y} + \mathbf{M}v_z \frac{\partial \mathbf{v}}{\partial z} = -\frac{1}{\rho}\frac{\partial p}{\partial x} + \nu \Delta \mathbf{v},$$

$$\text{div } \mathbf{v} = 0, \qquad (34)$$

$$\mathbf{y} = \mathbf{M}\mathbf{v}.$$

Variables $v_x$, $v_y$, $v_z$ are not signified as vectors because they are components of 3-dimensional vector $\mathbf{v} = (v_x, v_y, v_z) = \{\mathbf{v}_{ijk}\}$. But from the difference point of view they are the net functions. So $\mathbf{M}v_x$, $\mathbf{M}v_y$, $\mathbf{M}v_z$ are multiplications of matrix by vector.

Let a monotonizing operator $\mathbf{M}$ coincides (9). Note that in [17] the monotonizing operator was defined near the boundaries as well. So the results are different. The solution of the described problem at the 20×20×20 point mesh is shown at fig. 2. The solution $\mathbf{u}$ of the scheme of (8) type is nonmonotonic. It has 316 extremes in the interior of the cubic region. The solution $\mathbf{y}$ of the scheme of (10), (11) type is also nonmonotonic but it has only 112 extremes in the interior of the cubic region (an auxiliary function $\mathbf{v}$ has 312 extremes).

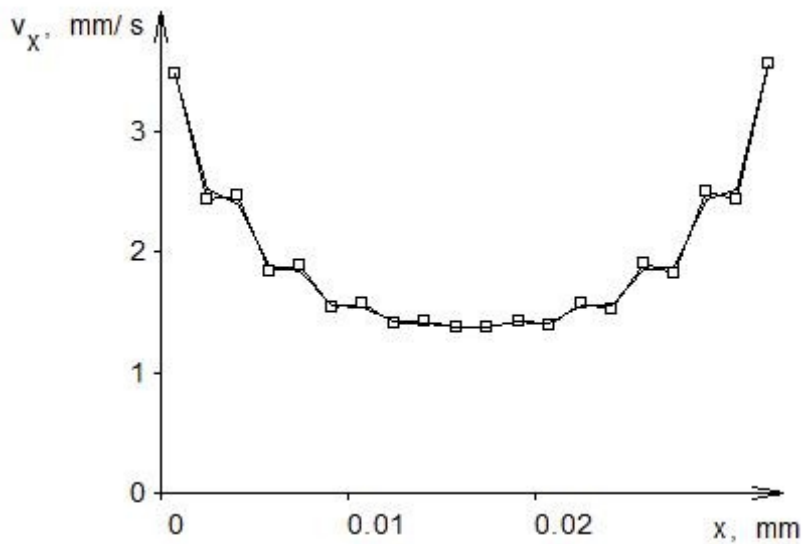

**Fig. 2** The solution of the problem of liquid flowing through the filter cell: the scheme (8) – markers, and the scheme (10), (11) – solid line. The velocity $v_x$ is shown on the line passing the centers of the holes. $L = 1/30$ mm, $\rho = 1$ mg/mm$^3$, $\nu = 1.002$ mm$^2$/s, $p_1$ – atmospheric

pressure, $p_0 = p_1 + 1000$ g/mm$^2$ s$^2$, $N=20$, the holes in the cube walls are the sum of boundary facets of the cells numbered $j,k = 6,7,...,15$

To characterize a sharpness of extremes let us introduce maximal change per step $a$ and maximal minimal change per step $b$ in a set S of extremes:

$$a = \max_{(i,j,k) \in S} \left( \max_{\alpha,\beta,\gamma \in H_{ijk}} |u_{ijk} - u_{\alpha\beta\gamma}| \right),$$
$$b = \max_{(i,j,k) \in S} \left( \min_{\alpha,\beta,\gamma \in H_{ijk}} |u_{ijk} - u_{\alpha\beta\gamma}| \right). \qquad (35)$$

Here $u$ is some net function and $H_{ijk}$ is a set of net points $\{(i+1,j,k), (i-1,j,k), (i,j+1,k), (i,j-1,k), (i,j,k+1), (i,j,k-1)\}$.

In the central part

$$S = \{(i,j,k): 6 \le i,j,k \le 15\} \qquad (36)$$

**u** has 48 extremes with maximal change per step $a=0.29$ in them and **v** has 4 extremes with maximal change per step $a=0.11$ in them. The maximal minimal change of **u** extremes in the central part is $b=0.0012$ and the same value for **v** is $b=1.4 \cdot 10^{-15}$.

So solutions of the monotonized schemes can be nonmonotonic. That is why we speak monotonization to be an improvement of monotonicity properties but don't speak monotonized schemes to be always monotonous.

**Conclusion.**

A previously proposed by the author method for improving of monotonicity properties of some classes of difference schemes is described in the article. The monotonized schemes can be formulated respectively to a similar function as the being monotonizing ones. Propositions concerning improvement of monotonicity properties are formulated and proved. Relation of such a monotonization with approximation is shown and the proposition concerning stability is proved. Many other properties are of more profound investigation which must be fulfilled more accurately. So it is out of the subject of this paper which reports some basic ideas of monotonized difference schemes constructing. Note that we speak not only about absolute monotonization but about improvement of monotonicity properties as well. I.e. the monotonized schemes generally speaking can be not monotonous in the sense of S.K. Godunov. But the propositions are new and they can be constructively applied in practical calculations.